\documentclass[11pt]{article}
\usepackage{amssymb, amsmath, amsthm}

\newcommand\blfootnote[1]{%
 \begingroup
 \renewcommand\thefootnote{}\footnote{#1}%
 \addtocounter{footnote}{-1}%
 \endgroup
}

\begin{document}
\parindent 0in
\parskip 1 em
\date{}
\title{\text{\bf{\normalsize{THE FUNCTIONAL DETERMINANT AND THE PARTITION FUNCTION}}}\\
 \text{\bf{\normalsize{IN GEOMETRIC FLOWS}}}}
\author{Christopher Lin}

\maketitle 

\begin{abstract}
 We propose the use of the functional determinant of geometric operators in constructing an entropy functional associated to geometric flows.  
 Our approach is based on the direct computation of the partition function, with a well-defined set of microstates and macrostates in the 
 canonical ensemble.  The approach is motivated by a fundamental enigma in Perelman's derivation of his famous $\mathcal{W}$-entropy.     
   The defining feature of our entropy is that the energy of each microstate in the partition function is invariant along the 
 associated geometric flow - a clue that could be inferred from Perelman's work.  Moreover, the monotonicity of our entropy 
 along the associated geometric flow is then a natural result in the statistical mechanics framework.      
 While we will not argue in a completely rigorous manner, we will use the formalism to derive an explicit formula for an 
 entropy associated to conformal flows on a closed surface based on the Polyakov formula for the determinant of the Laplacian.  We also 
 discuss possible extensions of our results to more general operators and manifolds.      
\end{abstract}

\blfootnote{\text{\it{Keywords and phrases}}\,: partition function, entropy, functional determinant, Laplacian, eigenvalues, Ricci flow, conformal 
deformation, Dirichlet energy \\
\text{\it{2010 Mathematics Subject Classification}}\,: 53A30, 53C44, 58J52, 82B30}

\section{Introduction}
\quad Functionals exhibiting monotonicity along solutions of partial differential equations have proven to be extremely useful.  
In a landmark paper \cite{P} that eventually contributed to the solution of the 
Poincar\`{e} conjecture, G. Perelman defined a functional he called the $\mathcal{W}$-entropy, which 
 is given by 
 \begin{equation}\label{W}
 \mathcal{W}(g,f,\tau) = \int_M \big(\tau(R+|\nabla f|^2)+ f - n\big) (4\pi\tau)^{-n/2}e^{-f}\,dv
\end{equation} 
on a closed manifold $M$ with any Riemannian metric $g$, any $f\in C^{\infty}(M)$, and any $\tau >0$.  
Perelman showed that along any solution to the system of differential equations
  \begin{equation}\label{system}
 \begin{cases}
  \frac{\partial g}{\partial t} = -2\text{Rc} \\
  \frac{\partial f}{\partial t} = -\Delta f + |\nabla f|^2 - R + \frac{n}{2\tau} \\
  \frac{d\tau}{dt} = -1\,,
 \end{cases}
\end{equation}
 which includes the Ricci flow $\frac{\partial g}{\partial t} = -2\text{Rc}$, the $\mathcal{W}$-entropy is non-decreasing.  More specifically, 
 he showed 
  \begin{equation}\label{mono}
 \frac{d\mathcal{W}}{dt} = \,2\tau\int_M \big|\text{Rc} \,+ \nabla^2 f - \frac{1}{2\tau}g\big|^2\,(4\pi\tau)^{-n/2}e^{-f}\,dv
\end{equation} 
along (\ref{system}).  

\quad The Lyapunov-type property of the $\mathcal{W}$-entropy above has far-reaching consequences in geometric analysis and in completing 
Hamilton's program of resolving the Poincar\`{e} conjecture.  However, a less talked-about, but nevertheless tentalizing fact is that 
  Perelman derived his $\mathcal{W}$-entropy from the formulas for the entropy in statistical mechanics.  To describe his derivation, let us 
  first recall some statistical physics below.
 
 \quad Recall that in statistical mechanics, the partition function for a thermodynamic system $\Gamma$ in a canonical ensemble is defined as 
\begin{equation}\label{par1}
 Z = \sum_i e^{-\beta E_i},
\end{equation}
where the summation is over all available ``microstates'' $\varphi_i$ of a system 
with $E_i$ the energy associated to $\varphi_i$, and $\beta = 1/\tau$ with $\tau$ as the temperature.  The partition function is the normalization 
factor in the \text{\it{Gibbs probability/measure}} 
\begin{equation}\label{prob}
 p_i = \frac{e^{-\beta E_i}}{Z}\,,
\end{equation}  
which is the probability that the system occupies the microstate $\varphi_i$ at thermal equilibrium (at temperature $\tau$).
In general the microstates may not be discrete, in which case we have the general formula 
\begin{equation}\label{genZ}
 Z = \int_\Omega e^{-\beta E(\omega)}\,d\omega \,,   
\end{equation}
where $\Omega$ is the space of microstates with $E:\Omega\longrightarrow \mathbb{R}$ the associated energy function, and 
$d\omega$ is a density of states measure on $\Omega$. Both the energy values of microstates and the density of 
states are independent of temperature.  The (equilibrium) entropy of $\Gamma$ at temperature 
$\tau$ is given by 
  \begin{equation}\label{eqgibbs}
 S = \log Z-\beta\frac{\partial}{\partial\beta}\log Z.
\end{equation}

\quad  In addition to microstates, there are also ``macrostates''.  
As a thermodynamic system, $\Gamma$ has 
a full set of macrostates describing its large-scale properties.  The temperature is a macrostate itself.     
A certain subset of macrostates are constraints that determine the energy value of each microstate, physical 
examples of macrostates of this type include the volume and the external magnetic field.   
From (\ref{eqgibbs}) one readily computes that 
\begin{equation}\label{dS}
 \frac{\partial S}{\partial\tau} = \frac{\sigma(\tau)^2}{\tau^3} \, ,
\end{equation}  
where $\sigma(\tau) \,= \sqrt{\,<(E - <E>)^2>}$ is the standard deviation of the energy with respect to the Gibbs measure (\ref{prob}).  Thus the 
entropy $S$ is non-decreasing when \text{\it{only}} the temperature is increased.

\quad In \cite{P}, Perelman treated $\tau$ as the temperature and it seems that he also conceived the Riemannian metrics $g$, smooth 
functions $f$ on an $n$-dimensional closed manifold as the energy-determining macrostates of some underlying thermodynamic system.  Then he declared      
\begin{equation}\label{lnZ}
 \log Z = \int_M (-f+\frac{n}{2})\,(4\pi\tau)^{-n/2}e^{-f}\,dv
\end{equation}    
without specifying what the underlying microstates and their energy values are (and hence also the partition function).  He proceeded to substitute 
(\ref{lnZ}) directly into (\ref{eqgibbs}), and differentiated along (\ref{system}) to arrive at the $\mathcal{W}$-entropy\footnote{his 
$\mathcal{W}$-entropy is negative $S$.} (\ref{W}).  Therefore Perelman did not 
propose the partition function $Z$, but instead its natural log.  This problem does not seem to generate much interest in the differential 
geometry circle,\footnote{the author is not sure about the situation in the theoretical physics community} as the author only knows of Xiang-Dong Li 
(see \cite{XDL}) who has addressed this issue directly.  One can of 
course choose to overlook it as pure coincidence, but the power that the $\mathcal{W}$-entropy 
has demonstrated in geometric analysis begs for a sound examination of its original motivation from statistical mechanics. 

\quad Motivated by the problem discussed above,  
we propose here a general scheme of defining partition functions associated 
to relevant geometric flows - based 
on the functional determinant of differential operators. This scheme has well-defined microstates as members of a functions space, 
with the energy of each microstate being the ``Dirichlet energy'' associated to the underlying operator.  The 
data that determines the operator (such as the Riemannian metric, additional smooth functions, etc.) then become the energy-determining 
macrostates.  Therefore, this is a scheme that follows the statistical mechanics prescription thoroughly.        

\quad  In view of formula (\ref{genZ}), a good reason for considering 
the functional determinant is the well-known 
formula 
 \begin{equation}\label{elementary}
 \int_{V}e^{-\langle\varphi,A\varphi\rangle}\,\mathcal{D}\varphi = (\det A)^{-1/2} 
\end{equation}
for a nonnegative-definite, self-adjoint operator $A$ on a finite $n$-dimensional inner product space $(V,\langle\,\,\,,\,\,\,\rangle$).  The measure 
$\mathcal{D}\varphi$ is the Euclidean measure $\pi^{-n/2}dx_1\cdots dx_n$, where 
$\langle\varphi,A\varphi\rangle = \sum_{i} \lambda_i x_i^2$ via the eigenfunction expansion 
$\varphi =\sum_i x_i\varphi_i$, $A\varphi_i = \lambda_i\varphi_i$.   To make mathematically rigorous sense of (\ref{elementary}) for an operator on an infinite dimensional function space, one needs regularization techniques 
such as the zeta-function regularization.  Physicists 
have long been using (\ref{elementary}) as a formal definition of the determinant of differential operators that they needed to compute 
in quantum field theory.  They have also explored its connection to partition functions via the zeta-function regularization (see the works \cite{GH} and \cite{H}, for example).\footnote{For a mathematical background on the zeta-function regularization, see the book \cite{BC} or the paper \cite{RS}.}
  Therefore, the idea of using the functional determinant to define partition functions is not at all new.  However, its use 
in defining entropy functionals in geometric evolution equations does not seem to have been explored before.

\quad The other key ingredient in our proposed scheme is the invariance of the Dirichlet energy $\langle\varphi,A\varphi\rangle$ of the 
differential operator $A$ along the associated geometric flow for each microstate $\varphi$.  This is the absolute essential property 
that we require between the operator $A$ and the associated geometric flow, and it is also the driving mathematical motivation 
behind our approach.  In general, we are seeking operators whose Dirichlet energy is invariant under a specified geometric flow, or vice-versa.  
This key ingredient stems from the observation that, in 
deriving the expression (\ref{W}) for the $\mathcal{W}$-entropy, Perelman 
took the \text{\it{total}} derivative of (\ref{lnZ}) with respect to $t$ (which is essentially $-\tau$) along the flow (\ref{system}).   
Since only the partial derivative is taken in (\ref{eqgibbs}), we can infer that if there really were a partition 
function (with well-defined microstates and energies) underlying the $\mathcal{W}$-entropy, then 
the system (\ref{system}) should be a flow along which the energy of each microstate $\varphi$ remains unchanged.  Moreover, 
this would mean that the monotonicity given by (\ref{mono}) is simply the natural result (\ref{dS}).

\quad  Despite the resemblence of (\ref{elementary}) to (\ref{genZ}), we cannot simply take (\ref{elementary}) as the definition of our partition function.  
The reason is because the measure $d\omega$ in (\ref{genZ}) must be independent of 
the macrostates,\footnote{besides just being independent of the temperature $\tau$} but (at least in the finite dimensional case) the measure $D\varphi$ in 
(\ref{elementary}) does depend on the operator $A$ and hence also on the associated energy-determining macrostates.  There is of course 
also the issue with the measure $D\varphi$ in (\ref{elementary}) not even being mathematically well-defined, but we will see how by forgoing this 
issue and pretend that the measure is there, we can still derive sensible formulas.  
Our approach is to somehow modify the measure that appears in (\ref{elementary}) into one that is macrostate-independent.    
We will do this by introducing a multiplicative factor to (\ref{elementary}) whose exact form 
is motivated by the invariance of the Dirichlet energy mentioned above.  The entropy is then calculated according to formula (\ref{eqgibbs}) by 
using the new partition function.

\quad At the present, we cannot verify that the $\mathcal{W}$-entropy indeed arises from the functional determinant approach described above.  In fact, 
in a later section we will explian how the $\mathcal{W}$-entropy cannot arise in \text{\it{exactly}}
 the way our scheme is carried out.  However, we will apply our 
scheme to some special casese.  First, we will work out the entropy 
formula for the Laplacian on closed surfaces, with the associated 
geometric flow being any conformal variation of the metric.    
This will not only yield a mathematically rigorous (at least within a conformal class) formula for the entropy, 
but it will also be explicitly computable.  The highlight of the end result is that in this case the entropy will be monotonic along any 
conformal variation of the metric in the temperature variable, hence corroborating (\ref{dS}) and is consistent with the idea of 
the invariance of Dirichlet energy along the flow.     
 
\quad Secondly, we will apply the same formalism and come up with the entropy formula for the so-called ``drifted Laplacian'', 
associated to any flow that contains a conformal deformation of the metric.  The additional data here is a smooth function, which 
determines the Dirichlet energy of the drifted Laplacian together with the metric.  
The entropy formula in this case is still more-or-less rigorous, but more work is needed to make it explicitly computable - at least on a 
closed surface.  We think that 
this entropy has potentially some connection to the $\mathcal{W}$-entropy since there both a metric and a function appears as determining 
data.        

\quad The organization of this note is as follows.  Section 2 contains work on the entropy for the Laplacian on closed surfaces - this section illustrates 
the methodology involved in the scheme described above, and should be considered to be the central focus of this note.  In Section 3, we attempt to 
generalize the results in Section 2 to the drifted Laplacian operator by reiterating the techniques illustrated in Section 2.  The discussions in 
Section 3 will actually be quite general, the drifted Laplacian is just one of the operators that we think is easier to tackle.  A brief, 
but noteworthy mention of conformal geometry is also at the end of Section 3.  In section 4 we 
return to our original motivation and discuss why our scheme needs to be modified in order that it even has a chance of giving rise to 
the $\mathcal{W}$-entropy.  Finally in section 5, we deem it necessary to give a physical interpretation of the objects appearing in our scheme.  
The author is not a physicist, so the reader is advised to be strongly spectical of the content in Section 5.  An appendix containing some 
calculations for the determinant in finite dimensions is attached at the end, with the hope of 
supplying credibility to some of the more sketchy arguments used in Sections 2 and 3.

\section{The case of the Laplacian on a Closed Surface}  
    \quad In this section, we will actually derive expressions for two related partition functions on a closed Riemann 
    surface - one for metrics in the same conformal 
class, and one for all metrics on the surface.  We will then compute entropy functions in the sense of (\ref{eqgibbs}) corresponding to each partition function.  
Up to a geometric constant, the partition function 
for all metrics generalizes that of metrics in the same conformal class.  The route we take here may seem redundant, but our goal is 
to demonstrate that different formulas turn out to be all consistent with each other.  This is particularly important due to the fact that we will 
be working with a measure on infinite dimensional function spaces that is not rigorously defined.

 \quad Incredibly, there are explicit mathematical results on the determinant of the Laplacian on a closed surface $M$.  For the remainder 
of this section, $A_g$ will denote the Laplacian $-\Delta$ acting on functions.  Note that $\Delta$ is defined as $\text{div}\,\nabla$ and thus 
nonpositive-definite.\footnote{In $\mathbb{R}^n$ with the flat metric, $\Delta$ is thus the sum of second partial derivatives in each coordinate.  Note 
that $\nabla$ is the gradient with respect to the metric.}  The following 
results appeared in various forms in the research literature, the versions below are extracted from the book \cite{BC}. 
\newtheorem{theorem}{Theorem}
\begin{theorem}\label{confvar}
 Suppose $g(t)$ is a smooth family of conformal metrics on a closed surface $M$, i.e. $\frac{\partial g}{\partial t} = \psi\,g$ for some function 
 $\psi(x,t)$.  Then 
 \[
  \frac{d}{dt}\log\det A_{g(t)} = -\int_M \psi\Big(\frac{R_{g(t)}}{24\pi}-\frac{1}{\text{Area}(g(t))}\Big)\,dA_{g(t)},
 \] 
 where $R_g$ is the scalar (twice the Gauss) curvature of $g$.  
\end{theorem}

\newtheorem{corollary}{Corollary}
\begin{corollary}\label{imp}
 Suppose two metrics $g$ and $h$ on a closed surface $M$ are related by $g = e^{\psi}h$.  Then 
 \[
  \log\det A_g - \log\det A_h = -\frac{1}{48\pi}\int_M\big(|\nabla\psi|_h^2+2\psi R_h\big)\,dA_h + \log\Big(\frac{\text{Area}(g)}{\text{Area}(h)}\Big).
 \]
\end{corollary}

\quad The formula in Corollary \ref{imp} was first 
derived by the Physicist A. Polyakov in \cite{Pv}, and such formulas have been called ``Polyakov formulas''.  In light of (\ref{genZ}), we define the quantity  
\begin{align}\label{Zg}
 Z_g(\beta) &= \int_\mathcal{H} e^{-\beta E_g(\varphi)}\,\mathcal{D}\varphi_g \notag\\
  &= (\det \beta A_g)^{-1/2},
\end{align}
where $\mathcal{H}$ is the domain of the Laplacian, $\beta = 1/\tau$, and 
$E_g(\varphi) = -\int_M \varphi\Delta\varphi\,dv = \int_M |\nabla\phi|^2 \,dv$ is the 
\text{\it{Dirichlet energy}} with all metric quantities dependent on the metric $g$.  
  The notation $\beta A_g$ is for the operator consisting of 
applying $A_g$ first and then mulitplying by $\beta$, although the order clearly does not matter.  
Since the space $C^{\infty}(M)$ of all smooth functions is dense in the domain of the Laplacian corresponding to any metric on $M$, we can 
assume $\mathcal{H} = C^{\infty}(M)$ for any metric $g$.  This last assumption means that our microstates are fixed with respect to changes in 
the energy-determining macrostates - a theme that will be repeated throughout.        

\quad  Note that the constant scaling of the metric by $\beta$ implies 
$\beta A_g = A_{\beta^{-1}g}$.  The first line in 
(\ref{Zg}) is a heuristic definition since we do not know how to construct the measure $\mathcal{D}\varphi_g$ rigorously, but the second 
line is well-defined mathematically by the zeta-function regularization.  We use the subscript $g$ in $Z_g(\beta)$ to denote 
that the measure $\mathcal{D}\varphi_g$ depends on the metric $g$, and as explained in the introduction it means that 
(\ref{Zg}) should not be considered as a true partition function.  The rationale for the dependence of $\mathcal{D}\varphi_g$ on the metric 
is that this is true for the finite dimensional case (see Appendix).  We will modify the measure $\mathcal{D}\varphi_g$ so to that it 
eventually takes a form that is metric-invariant.     

\quad  We want to point out a confusion that may arise as a result of the notations used here, and it reinforces the importance of the dependence of the measure 
$\mathcal{D}\varphi_g$ on $g$.  As operators, $\beta A_g$ and $A_{\beta^{-1}g}$ are equal, and we are computing the Dirichlet energy of $A_{\beta^{-1}g}$ with respect to the volume form $dv_g$ in (\ref{Zg}).  We are \text{\it{not}} computing the Dirichlet energy of $A_{\beta^{-1}g}$ in (\ref{Zg}) 
with respect to the volume 
form $dv_{\beta^{-1}g}$, in which case it would be equal\footnote{This is 
a special case of the conformal invariance of the Dirichlet energy of the Laplacian on surfaces, as we will see below.} to the Dirichlet energy of $A_g$ with respect to the volume form $dv_g$.  To put things in perspective, from definition (\ref{Zg}) we have 
\begin{align}
 (\det \beta A_g)^{-1/2} &= \int_{\mathcal{H}} e^{-\beta\langle\varphi, A_g \varphi\rangle_g}\,\mathcal{D}\varphi_g \notag\\
            &= \int_{\mathcal{H}} e^{-\langle\varphi, A_{\beta^{-1}g}\varphi\rangle_g}\,\mathcal{D}\varphi_g \notag
\end{align}
On the other hand, we must also have 
\begin{align}
 (\det \beta A_g)^{-1/2} &= (\det A_{\beta^{-1}g})^{-1/2} \notag \\
                         &= \int_{\mathcal{H}} e^{-\langle\varphi, A_{\beta^{-1}g} \varphi\rangle_{\beta^{-1}g}}\,\mathcal{D}\varphi_{\beta^{-1}g} \notag\\
                         &= \int_{\mathcal{H}} e^{-\langle\varphi, A_g \varphi\rangle_g}\,\mathcal{D}\varphi_{\beta^{-1}g}. \notag
\end{align}
For consistency's sake, we must then have 
\begin{equation}\label{classic}
 \int_{\mathcal{H}} e^{-\langle\varphi, A_{\beta^{-1}g}\varphi\rangle_g}\,\mathcal{D}\varphi_g = 
 \int_{\mathcal{H}} e^{-\langle\varphi, A_g \varphi\rangle_g}\,\mathcal{D}\varphi_{\beta^{-1}g},
\end{equation}
which can be easily verified in the finite-dimensional case (see Appendix).  

\quad First we consider metrics in the same conformal class on a closed surface.  
 Suppose $g(t)$ is a smooth family of metrics in the same conformal class on a closed surface $M$.  The vital observation here is that 
 the Dirichlet energy is invariant under conformal deformations.  To see this, suppose $\frac{\partial g}{\partial t} = \psi \,g$ 
is a conformal deformation on $M$.  For a 
fixed function $\varphi$, we have the elementary formula 
\begin{equation}\label{varnorm}
 \frac{\partial}{\partial t}|\nabla\varphi|^2_g = -\frac{\partial g}{\partial t}(\nabla\varphi,\nabla\varphi)
\end{equation}    
whose proof we will skip.  Using (\ref{varnorm}) and another well-known formula 
$\frac{\partial}{\partial t}dv = \frac{1}{2}\text{tr}_g\big(\frac{\partial g}{\partial t}\big)$, we easily compute that\footnote{Note that this, and 
a more general variational formula to appear in the next section, is not computing the variation of \text{\it{eigenvalues}} of the operator.  Relevant, 
and quite interesting work on the variation of eigenvalues along geometric flows can be found in \cite{C}, \cite{Li}, and \cite{Ma}.  We think however, 
that this latter approach may figure in perahps a different approach to establishing the $\mathcal{W}$-entropy.} 
\[
 \frac{d}{dt}\int_M |\nabla\varphi|^2\,dv = -\int_M \psi|\nabla\varphi|^2\,dv + \int_M |\nabla\varphi|^2\,\psi\,dv = 0.
\]
We will take advantage of this fact in the following argument.  
Since $E_g$ is invariant under conformal deformations of $g$ and 
therefore the only object dependent on $g$ in the left-hand side of (\ref{Zg}) is the measure $\mathcal{D}\varphi_g$, we must have
\begin{equation}\label{must}
 \frac{\partial Z_{g(t)}(\beta)}{\partial t} = \int_{\mathcal{H}} e^{-\beta E_g(\varphi)}\,\frac{\partial}{\partial t}\mathcal{D}\varphi_g.
\end{equation} 
Then by the identification (\ref{Zg}) and using Theorem \ref{confvar}, we compute that    
\begin{align}\label{disclosed}
 \frac{\partial}{\partial t}\log Z_{g(t)}(\beta) 
  &= -\frac{1}{2}\frac{\partial}{\partial t}\log\det(\beta A_{g(t)}) \notag\\
  &= -\frac{1}{2}\Big(-\int_M \psi\Big(\frac{R_{\tilde{g}(t)}}{24\pi}-\frac{1}{\text{Area}(\tilde{g}(t))}\Big)dv_{\tilde{g}(t)}\Big) \notag\\ 
  &= \frac{1}{2}\int_M \psi\Big(\frac{R_{g(t)}}{24\pi}-\frac{1}{\text{Area}(g(t))}\Big)dv_{g(t)} \notag\\
  &= -\frac{1}{2}\frac{d}{d t}\log\det A_{g(t)} \notag\\
  &= \frac{d}{d t}\log Z_{g(t)}(1),
\end{align}
where we defined $\tilde{g}(t) = \beta^{-1}g(t)$ which still satisfies $\frac{\partial\tilde{g}}{\partial t} = \psi\tilde{g}$, and note that the constant 
factor $\beta^{-1}$ scales out of the curvature and volume so as to cancel out with the scaling of the volume form.  
     
\quad Let us define 
\begin{equation}\label{desired}
 Z(\beta) = Z_g(\beta)/Z_g(1) 
\end{equation} 
and further define the new measure 
\begin{equation}\label{newm}
 \mathcal{D}\varphi = \frac{\mathcal{D}\varphi_g}{Z_g(1)},
\end{equation}
so that we can formally write 
\begin{equation}\label{formal}
 Z = \int_{\mathcal{H}}e^{-\beta E_g(\varphi)}\,\mathcal{D}\varphi.
\end{equation}    

\quad  Let us write $\Phi(t) = \frac{d}{d t}\log Z_{g(t)}(1)$, then it follows from (\ref{must}) and (\ref{disclosed}) that 
\begin{align}\label{rearrange}
 \int_{\mathcal{H}} e^{-\beta E_g(\varphi)}\,\frac{\partial}{\partial t}\mathcal{D}\varphi_g &= \Phi(t)Z_{g(t)}(\beta)\notag\\
  &= \Phi(t)\int_{\mathcal{H}} e^{-\beta E_g(\varphi)}\,\mathcal{D}\varphi_g.
\end{align}  
We can go a step further and postulate from (\ref{rearrange}) that 
\begin{equation}\label{soas}
 \frac{\partial}{\partial t}\mathcal{D}\varphi_g = \Phi(t)\mathcal{D}\varphi_g
\end{equation} 
under any conformal deformation of $g$.\footnote{Note however, that by our identification (\ref{Zg}) which resulted directly in (\ref{disclosed}), 
equation (\ref{rearrange}) is valid 
without assuming (\ref{soas}).}  Then by differentiating (\ref{newm}) directly and using (\ref{soas}), we have  
\begin{align}
 \frac{\partial \mathcal{D}\varphi}{\partial t} &= \frac{\partial}{\partial t}\frac{\mathcal{D}\varphi_g}{Z_g(1)} \notag\\
  &= -\frac{1}{Z_g(1)^2}\frac{d Z_g(1)}{d t}\mathcal{D}\varphi_g + 
  \frac{1}{Z_g(1)}\frac{\partial\mathcal{D}\varphi_g}{\partial t} \notag\\
  &= \frac{1}{Z_g(1)}\Big(-\frac{d}{d t}\log Z_g(1)\,\mathcal{D}\varphi_g + \Phi(t)\mathcal{D}\varphi_g\Big) = 0.  
\end{align}
Thus the measure (\ref{newm}) is invariant under conformal deformations of the metric as we had 
hoped to find, assuming that (\ref{soas}) holds.\footnote{Note that a measure like (\ref{newm}) 
can always be defined, but its conformal invariance is rooted in the result (\ref{disclosed}).}  With this, $Z$ as rigorously 
defined by (\ref{desired}) and formally given by (\ref{formal}) can now be considered as a true partition function.    

\quad Let us proceed to compute the entropy associated to the partition function defined by (\ref{desired}).  By Corollary \ref{imp} and using 
the Gauss-Bonnet formula, we compute that 
\begin{equation}\label{mar}
 \log Z = \Big(\frac{1}{2} - \frac{\chi(M)}{12}\Big)\log \beta.
\end{equation}
Then it follows readily that 
\begin{equation}\label{turnsout}
 S = \Big(\frac{1}{2} - \frac{\chi(M)}{12}\Big)\big(\log\beta - 1\big).  
\end{equation}
Therefore the partition function and its entropy only depend on the topology of the surface $M$.  
This is a bit surprising in light of the fact that we only have the respective conformal invariance of the Dirichlet energy $E_g$ 
and the measure $\mathcal{D}\varphi$.  
The correct interpretation here is that although the partition function $Z$ (and hence the corresponding entropy $S$) defined above is 
a-priori dependent on the conformal class of the metric, as the metric moves across conformal classes the  
effect of the conformal variation in the Dirichlet energy $E_g$ and that of the measure $\mathcal{D}\varphi$ cancel eachother.
  
\quad From (\ref{turnsout}) we have  
\begin{equation}\label{varS}
 \frac{\partial S}{\partial \tau} = \Big(\frac{\chi(M)}{12} - \frac{1}{2}\Big)\frac{1}{\tau},
\end{equation}
which has a definite sign for all $\tau >0$ except at the critical value of $\chi(M) = 6$ where it equals to zero.  
In particular, when $M$ is in addition orientable we have 
\begin{equation}\label{varS2}
 \frac{\partial S}{\partial \tau} = -\Big(\frac{2+\,\text{genus}}{6}\,\Big)\frac{1}{\tau} <0.
\end{equation}
Thus for a closed Riemann surface, the entropy (\ref{turnsout}) is definitely monotonic in the 
temperature $\tau = 1/\beta$.  It is also 
important to point out that from our construction, (\ref{varS}) is also attained  
by taking the total variation in $\tau$ along conformal variations of the metric, i.e. along a flow of the form 
\begin{equation}\label{latenight}
 \begin{cases}
  \frac{\partial g}{\partial t} = \psi\,g \\
  \frac{d \tau}{d t} = -1
 \end{cases}
\end{equation}
analogous to (\ref{system}) for the $\mathcal{W}$-entropy.  We would like to point out that the Ricci flow is also a conformal flow in dimension $2$, as it 
is given by $\frac{\partial g}{\partial t} = -R\,g$, where $R$ is the scalar curvature of $g$.     

\quad The work above is just a warm-up to what we are really after, since the entropy (\ref{turnsout}) is for a thermodynamic system 
whose energy-determining macrostates (metrics in a fixed conformal class) all give the same microstate energy - quite a 
dull system indeed.  We want to enlarge the space of energy-determining 
macrostates to all metrics on a closed surface.  To this end, we look for a more general expression for the entropy in the form of (\ref{formal}) that contains 
a metric-independent (not just conformally-invariant) density measure $\mathcal{D}\varphi$.  
The following argument generalizes the one that lead to the previous entropy.  Based on its validity in the finite dimensional case, we will assume that 
\begin{equation}\label{justlike}
 \frac{\partial}{\partial t}\mathcal{D}\varphi_g = \mathcal{F}\,\mathcal{D}\varphi_g
\end{equation}            
under any variation of the metric, for a function $\mathcal{F}$ independent of $\varphi\in\mathcal{H}$.  Then we 
consider the general ansatz for the desired partition function as $Z = \eta(g) \, Z_g(\beta)$ for some function $\eta$ that depends only on the metric $g$, 
and where $Z_g(\beta)$ is as in (\ref{Zg}).  
In terms of integration on $\mathcal{H}$, we can then write 
\[
 Z = \int_{\mathcal{H}} e^{-\beta E_g(\varphi)}\,\eta(g)\,\mathcal{D}\varphi_g.
\]
Thus we have a new measure $\eta(g)\,\mathcal{D}\varphi_g$, and the requirement that it be metric-independnt means that 
\begin{equation}\label{defeqn}
 \frac{d\eta}{dt} + \mathcal{F}\,\eta = 0
\end{equation}       
for all variations $g(t)$ of the metric.  To see what $\mathcal{F}$ has to be, we compute
\begin{align}
 \frac{\partial}{\partial t}\log Z_{g(t)}(\beta) &= \Big(-\beta\int_{\mathcal{H}} \frac{\partial E_g(\varphi)}{\partial t}\,e^{-\beta E_g(\varphi)}\,\mathcal{D}\varphi_g 
 \,+\, \int_\mathcal{H} e^{-\beta E_g(\varphi)} \,\mathcal{F}\,\mathcal{D}\varphi_g\Big)/Z_{g(t)}(\beta) \notag\\
 &= -\beta\Big< \frac{\partial E_g(\varphi)}{\partial t} \Big>_{\beta} + \mathcal{F}, \notag 
\end{align}
where $\big<\,\,\,\,\big>_\beta$ denotes the 
average (at temperature $\tau = \beta^{-1}$) with respect to the Gibbs measure, which can be used since 
the $\eta$ factors can be cancelled out.  Then we see that 
\begin{equation} \label{iffy}
 \mathcal{F} = \beta\Big<\frac{\partial E_g(\varphi)}{\partial t}\Big>_{\beta} + \frac{\partial}{\partial t}\log Z_{g(t)}(\beta).
\end{equation}
Since $\mathcal{F}$ must also be independent of $\beta$, (\ref{iffy}) must hold for all $\beta$.  Setting $\beta=1$ we then have 
\begin{equation}\label{better}
 \mathcal{F} = \Big<\frac{\partial E_g(\varphi)}{\partial t}\Big>_{1} + \Phi(t),
\end{equation} 
where $\Phi(t)$ is as in (\ref{rearrange}).  Note that when the variation in the metric is conformal, 
$\frac{\partial E_g(\varphi)}{\partial t}=0$ and we immediately recover (\ref{soas}).

\quad In order to incorporate (\ref{better}) in solving for $\eta$ in (\ref{defeqn}), we need $\Big<\frac{\partial E_g(\varphi)}{\partial t}\Big>_{1}$ 
to be the variation of some function of metrics on the space of metrics.\footnote{This can be seen as an integrability requirement on the space of metrics.}  
However, it is difficult to show this in general.  Instead, we seek to find another way to express $\mathcal{F}$.  To this end, we will now make an 
additional assumption that between any two metrics $g_0$ and $g$, we have 
\begin{equation}\label{additional}
 \mathcal{D}\varphi_g = J(g_0, g)\mathcal{D}_{g_0}
\end{equation}   
for some positive function $J$ dependent on $g_0$ and $g$ but independent of $\varphi\in\mathcal{H}$.  The function $J$ is reminiscent of the Jacobian in 
calculus (see the Appendix for more on this).  With a fixed metric $g_0$ chosen, $J(g_0, g)$ is a function of the metric $g$.  Then by varying 
the metric, we compute that 
\begin{align}
 \frac{\partial}{\partial t}\mathcal{D}\varphi_{g(t)} &= \frac{\partial J}{\partial t}\,\mathcal{D}\varphi_{g_0}\notag\\
  &= \frac{\partial \log J }{\partial t}\, \mathcal{D}\varphi_{g(t)}\notag
\end{align}     
in view of (\ref{additional}).  Note that (\ref{additional}) also directly implies (\ref{justlike}), out of which we can now conclude
\begin{equation}\label{hao}
 \mathcal{F} = \frac{\partial \log J }{\partial t}.
\end{equation} 
Therefore, we have succeeded in identifying a function $J$ of the metric whose variation in the space of metrics gives us $\mathcal{F}$.  Then 
it is clear that $\eta(g) = J(g_0,g)^{-1}$ solves (\ref{defeqn}) for all variations $g(t)$ of the metric.  Note that by (\ref{additional}), 
\begin{equation}\label{transition}
 J(g_0,g) = J(g_0,g_1)J(g_1,g)
\end{equation}
for any three metrics $g_0, g_1$, and $g$.  Although $J(g_0,g)$ depends on the base 
metric $g_0$ chosen, (\ref{transition}) implies that (\ref{hao}) is independent of the base metric - as it should be.  Therefore, the more general 
entropy can be written as  
\begin{equation*}
 Z = \frac{Z_{g}(\beta)}{J(g_0,g)}
\end{equation*}     
with $g_0$ some fixed metric, which generalizes (\ref{desired}).  In fact, by (\ref{additional}) we see that 
\[
 Z_{g}(1) = \int_{\mathcal{H}}e^{-E_g(\varphi)}\,J(g_0,g)\mathcal{D}\varphi_{g_0}
\] 
and it allows us to rewrite the general partition function as 
\begin{equation}\label{final}
 Z = \frac{Z_{g}(\beta)}{Z_{g}(1)}\int_\mathcal{H}e^{-E_g(\varphi)}\,\mathcal{D}\varphi_{g_0}.
\end{equation}
The partition function (\ref{final}) is our final desired product, which involves a measure 
\[
 \mathcal{D}\varphi = \frac{\int_\mathcal{H}e^{-E_g(\varphi)}\,\mathcal{D}\varphi_{g_0}}{Z_{g}(1)}\mathcal{D}\varphi_g
\]  
that is invariant under the choice of macrostate $g$.  We speculate that a perturbative technique might be available 
to calculate this general partition function, since one clearly sees 
from (\ref{final}) the need to expand $E_g(\varphi)$ about $g_0$.  

\quad On a closed surface with metric $g$, by taking the partial derivative with respect to 
$\beta$ and using (\ref{turnsout}) the corresponding entropy to (\ref{final}) is then  
\begin{equation}\label{fentro}
 S =  \Big(\frac{1}{2} - \frac{\chi(M)}{12}\Big)\big(\log\beta - 1\big) + \log\int_{\mathcal{H}}e^{-E_g(\varphi)}\,\mathcal{D}\varphi_{g_0},
\end{equation}     
 and we see the presence of the extra geometric term.\footnote{Note that this is a 
\text{\it{relative}} geometric term that depends on a fixed metric $g_0$, hence (\ref{final}) should be seen as a \text{\it{relative entropy}}.}  
Note that via the conformal invariance of the Dirichlet energy $E_g(\varphi)$, the extra geometric term must be constant along conformal 
variations of the metric.  Thus we also arrive at the 
same formula (\ref{fentro}) by taking the total derivative with respect to $\beta$ along the flow (\ref{latenight}).  
  
\quad Furthermore, by taking the partial derivative with respect to $\tau$ we still have 
(\ref{varS}) and (\ref{varS2}), which also hold by taking the total derivative with 
respect to $\tau$ along the flow (\ref{latenight}).  In particular, for a closed Riemann surface 
the entropy (\ref{fentro}) is still monotonic along (\ref{latenight}).        
The entropy (\ref{fentro}) again reflects the key property that certain geometric flows (any conformal flow in this case) preserve the energy of 
each microstate of the underlying system.  We also want to point out that even if the entropy (\ref{fentro}) cannot be made mathematically 
rigorous due to the presence of the extra geometric term, its variation in temperature (\ref{varS}) is completely well-defined.     

\quad A few words about the computability of (\ref{fentro}).  Note that at the same temperature, 
the entropy (\ref{fentro}) is constant for metrics $g$ in the same conformal class, with any base metric $g_0$.  
Moreover, if $g$ is in the same conformal class as $g_0$, then the entropy becomes
\begin{equation}\label{sameclass}
  S = \Big(\frac{1}{2} - \frac{\chi(M)}{12}\Big)\big(\log\beta - 1\big) -\frac{1}{2}\log\det A_{g_0}
\end{equation}
by the invariance of the Dirichlet energy.  Although we cannot compute (\ref{sameclass}) explicity, if we fix another metric $g_1 = e^{\psi}g_0$ for some 
function $\psi$, then by Corollary \ref{imp} we see that 
\[
 S_{g_1} - S_{g_0} = \frac{1}{96\pi}\int_M\big(|\nabla\psi|_h^2+2\psi R_{g_0}\big)\,dA_{g_0} - 
                       \frac{1}{2}\log\Big(\frac{\text{Area}(g_1)}{\text{Area}(g_0)}\Big).
\]  

\quad One last remark is in order before we end this section.  We want to 
stress again that the monotonicity of our entropy (\ref{fentro}) exhibited by (\ref{varS}) should 
be seen as a testamant to the validity of our scheme, in spite of the lack of mathematical rigor in our arguments.  It is also 
worth pointing out that if we had not modified $Z_g(\beta)$ by a multiplicative factor and just took the functional determinant 
in (\ref{Zg}) as the definition of our partition function, then the resulting entropy $S$ according to (\ref{eqgibbs}) will not 
be monotonic along the flow (\ref{latenight}).  An obvious issue here is that the signs do not seem 
to match up, as the sign in (\ref{dS}) is always nonnegative while typically\footnote{for example, closed Riemann surfaces} the sign in (\ref{varS}) is 
nonpositive.  This issue is fundamentally rooted in the fact\footnote{It is a simple exercise to check this using Theorem \ref{confvar}.}
 that if we scale a fixed metric $g$ by $\tilde{g}(\tau) = \tau g$, then
 \[
 \frac{d}{d\tau}\log\det A                 
  = -\frac{1}{\tau}\Big(\frac{\chi(M)}{6} -1\Big),
 \]
 where $A = -\Delta_{\tilde{g}(\tau)}$ and again we see the critical value of $\chi(M)=6$.  In particular, $\frac{d}{d\tau}\log\det A = (g+2)/3\tau >0$ when 
 $M$ is a closed Riemann surface.     
 This is counter-intuitive, as scaling a metric larger will lower the eigenvalues 
of the Laplacian and hence should also decrease the value of the determinant, since it is the product of eigenvalues.  Somehow the zeta-function 
regularization has destroyed this intuition, but nevertheless the monotonicity of the entropy is still preserved.

\section{Other Operators} 
\quad The more ambitious task is to follow the same scheme as in the last section for obtaining the expression for entropy on higher dimensional manifolds with 
possibly different operators.  More precisely, similar to (\ref{formal}) we want a partition function in the form 
\[
 Z = \int_{\mathcal{H}}e^{-\beta E(\varphi)}\,\mathcal{D}\varphi,
\] 
where $E(\varphi) = \langle\varphi,A\varphi\rangle$ is the Dirichlet energy and $\mathcal{D}\varphi$ is a fixed measure 
on the domain $\mathcal{H}$ of a more general operator $A$.  A natural operator to consider would be the drifted 
Laplacian $\Delta_f = \Delta - \langle\nabla f,\nabla\,\,\rangle$ on a Riemannian manifold $M$ with metric $g = \langle\,\,\,,\,\,\,\rangle$ and 
``potential function'' $f\in C^{\infty}(M)$.

\quad As in the last section, we wish to first identify invariance properties of the ``Dirichlet energy'' for such operators.  
The following generalizes the conformal invariance of the Dirichlet energy of the Laplacian on a closed surface. 
 Let $M$ be a closed manifold equipped with a Riemannian metric $g$, and two functions $f$, $V\in C^{\infty}(M)$.  We consider in general the 
drifted Schr\"{o}dinger operator $A = -\Delta + \langle\nabla f,\,\nabla\,\,\rangle + V$.  Suppose we have a system of evolution equations
\begin{equation}
 \begin{cases}
  \frac{\partial g}{\partial t} = h(g,f,V)\\
  \frac{\partial f}{\partial t} = F_1(g, f,V) \\
  \frac{\partial V}{\partial t} = F_2(g, f,V)
 \end{cases}
\end{equation}   
Then we compute 
\begin{align}\label{realization}
 \frac{d}{dt}\int_M \phi A\phi \,du &= \frac{d}{dt} \int_M \big(|\nabla\phi|^2 + V\phi^2\big)\,du \notag\\
  &= \int_M \big(-h(\nabla\phi,\nabla\phi)+ F_2\,\phi^2\big)\,du + 
  \int_M \big(|\nabla\phi|^2 + V\phi^2\big)\big(\frac{1}{2}\text{tr}\,_g(h) - F_1\big)\,du.
\end{align}
Consider a conformal deformation of the metric $g$, i.e. $h = \psi\,g$ for some function $\psi(x,t)$.  Then 
$\text{tr}\,_g(h) = n\psi$, where $n$ is the dimension of $M$.  In view of (\ref{realization}), we see that by letting 
$F_1 = \big(\frac{n}{2} - 1\big)\psi$ and $F_2 = -\psi\, V$ we can make the derivative of the Dirichlet energy vanish.  
Thus we have the following result.

\newtheorem{proposition}{Proposition} 
\begin{proposition}\label{main}
 On a compact manifold $M$ of dimension $n$, if $(g,f,V)$ satisfies a system of evolution equations of the form
 \begin{equation}\label{bestsystem}
  \begin{cases}
   \frac{\partial g}{\partial t} = \psi\,g \\
   \frac{\partial f}{\partial t} = \big(\frac{n}{2} -1\big)\psi \\
   \frac{\partial V}{\partial t} = -\psi\, V
  \end{cases}
 \end{equation}
 for some function $\psi(x,t)$, then the Dirichlet energy of the associated drifted Schr\"{o}dinger operator $A$ is invariant under the system.  
\end{proposition}  
We would like to point out that when $n\geq 3$,  $\frac{\partial g}{\partial t} = -R\,g$ is the Yamabe flow.\footnote{see \cite{BC} for a background}   

\quad For 
simplicity we will just investigate the case of the drifted Laplacian operator $A = -\Delta_f$.  
Let us denote the data $(g,f)$ by $\theta$, so that we can more precisely denote $A = A_{\theta}$.  Moreover, 
we will only discuss the analogous formula to (\ref{fentro}), since what we really want is an invariant measure $\mathcal{D}\varphi$ within 
some class of $\theta$, with $E_{\theta}(\varphi)$ nonconstant within the class.  

\quad Again we assume that 
$\mathcal{D}\varphi_\theta = J(\theta_0,\theta)\mathcal{D}\varphi_{\theta_0}$
for some function $J(\theta_0,\theta)$ independent of $\varphi$, for any $\theta = (g,f)$ and a fixed $\theta_0 = (g_0,f_0)$.  Then by the same 
argument as in the last section, we readily arrive at the general partition function 
\begin{equation}\label{further}
 Z = \frac{Z_{\theta}(\beta)}{Z_\theta (1)}\int_{\mathcal{H}} e^{-E_{\theta}(\varphi)}\,\mathcal{D}\varphi_{\theta_0}
\end{equation}  
with a fixed macrostate $\theta_0$, and where 
$Z_\theta (\beta)$ is defined by the functional determinant of $\beta A_\theta$ just like 
in (\ref{Zg}).\footnote{The zeta-function regularized determinant actually works for any self-adjoint elliptic operator on a closed manifold of 
any dimension.  The drifted Laplacian is such an operator.}  
The partition function (\ref{further}) holds for $\theta = (g,f)$ with any metric $g$ and any smooth function $f$ on $M$.  
Therefore just as in the last section, the underlying (though not rigorous) measure 
in the canonical ensemble is 
\[
 \mathcal{D}\varphi = \frac{\int_{\mathcal{H}} e^{-E_{\theta}(\varphi)}\,\mathcal{D}\varphi_{\theta_0}}{Z_\theta (1)}\mathcal{D}\varphi_\theta,
\]
which is independent of $\theta=(g,f)$ with any metric $g$ and any smooth function $f$ on $M$.
  
\quad By Proposition \ref{main} we see that any variation in $\theta$ consisting of 
the first two equations of system (\ref{bestsystem}) leaves the energy $E_\theta(\varphi)$ of each microstate $\varphi$ 
unchanged.  Therefore we can consider the flow 
\begin{equation}\label{superman}
 \begin{cases}
  \frac{\partial g}{\partial t} = \psi\,g \\
   \frac{\partial f}{\partial t} = \big(\frac{n}{2} -1\big)\psi \\
   \frac{d\tau}{dt} = -1,
 \end{cases}
\end{equation} 
which is analogous to (\ref{system}).  Then taking the total variation of (\ref{further}) 
in the variable $\beta$ along (\ref{superman}), the entropy according to (\ref{eqgibbs}) 
becomes
\begin{equation}\label{newS}
 S = -\frac{1}{2}\Big(\log\frac{\det\beta A_\theta}{\det A_\theta} - \beta\frac{d}{d\beta}\big(\log\frac{\det\beta A_\theta}{\det A_\theta}\big) \Big) 
 + \log\int_{\mathcal{H}}e^{-E_{\theta}(\varphi)}\,\mathcal{D}\varphi_{\theta_0}.
\end{equation}
The entropy formula (\ref{newS}) holds for the drift Laplacian $A_\theta$ on any closed $n$-dimensional manifold.  At the present we 
cannot compute (\ref{newS}) further, because to our knowledge 
the Polyakov formula $\log\frac{\det\beta A_\theta}{\det A_\theta}$ for the drifted Laplacian is not known - even on closed surfaces.  
However, we believe that such a Polyakov formula must exist for closed surfaces just as in the case of the Laplacian.  If this 
Polyakov formula is computed, then we can also use it to prove that the entropy (\ref{newS}) is monotonic along 
the flow (\ref{superman}).\footnote{Note that by the invariance of the energy $E_\theta(\varphi)$, the last term in (\ref{newS}) also 
vanishes along (\ref{superman}).  Therefore, the variation of the entropy again dependents only  on the Polyakov formula of the operator 
in question, just as in the previous case of the Laplacian.}          

\quad Let us elaborate on the conjecture mentioned above for the drifted Laplacian on closed surfaces.  It is elementary to check that 
\begin{align}
\beta A_{\theta} &= -\beta\Delta_{g} + \beta\langle\nabla^{g}f,\nabla^g\,\,\,\rangle_{g} \notag\\
                 &= -\Delta_{\beta^{-1}g} + \langle\nabla^{\beta^{-1}g}f,\nabla^{\beta^{-1}g}\hskip 0.4cm\rangle_{\beta^{-1}g} = A_{(\beta^{-1}g,f)}
\end{align}
since the gradient scales as $\nabla^{cg} = c^{-1}\nabla^g$ for any $c>0$.  Then taking the partial derivative(s) with respect to $\beta$ of 
$\log\det \beta A_{\theta}$ amounts to varying the metric $\beta^{-1}g$ conformally in $\beta$ 
with $g$ fixed (and $f$ fixed as well) just as in the case of the ordinary Laplacian.  Thus we expect to get  
formulas similar to those in Theorem \ref{confvar} and Corollary \ref{imp}, and since the fixed potential function $f$ is 
outstanding throughout, we expect $f$ to also appear explicitly in the final formulas for $S$ and $\frac{\partial S}{\partial \tau}$.  
This suggests 
that if explicitly worked out, 
$\log Z_{\theta} = (-1/2)\log\det A_{\theta}$ could resemble\footnote{but not equal to.  See next section for more details.} Perelman's formula (\ref{lnZ}).    

\quad On the other hand, instead of considering all metrics on a closed surface, it also makes sense to consider the set of energy-determining 
macrostates as $[g]\times C^{\infty}(M)$, with $[g]$ a fixed conformal class.      
This means that the Dirichlet energy $E_{\theta} (\varphi)$ will only vary nontrivially with 
the potential function $f$.  
In fact, by deforming the metric inside $[g]$ and simultaneously deforming the potential function $f$ arbitrarily, 
(\ref{realization}) implies 
\begin{equation}\label{skywalker}
 \frac{\partial E_{\theta(t)}(\varphi)}{\partial t} = -\int_M \frac{\partial f}{\partial t} |\nabla\varphi|^2\,du. 
\end{equation}
Therefore we think that this is the setting in which the entropy (\ref{newS}) is most amenable to computation, 
since (\ref{skywalker}) is a manageable first-order approximation 
to the Dirichlet energy.  In any case, the precise form of formula (\ref{newS}) will hinge on the Polyakov formulas for the drifted Laplacian.

\quad The success we have had with conformal deformations also points to another general direction.  There is a Polyakov formula for 
any second-order, conformally covariant operator in dimension $4$.  The formula was discovered by Branson and \O rsted, and like the Polyakov 
formula in Corollary \ref{imp} it is an integral expression containing, now, new curvature terms.  
We refer to the lecture notes 
\cite{AC} of A. Chang for an extensive coverage on this subject.\footnote{The short, but very informative review \cite{RM} by R. Mazzeo is a 
helpful start.}  It will be a good idea to apply our scheme 
to such conformally covairant operators and derive the corresponding entropy functionals associated to conformal deformations of the metric.  This 
potentially fruitful direction may shed some new light on the field of conformal geometry by linking it to statistical mechanics.

\section{Further Speculations on the $\mathcal{W}$-entropy}
\quad The $\mathcal{W}$-entropy has a salient feature when (\ref{mono}) vanishes, which occurs along a smooth 
one-parameter family of $(g,f,\tau)$ satisfying the equation
\begin{equation}\label{soliton}
 \text{Rc} \,+ \nabla^2 f - \frac{1}{2\tau}g = 0.
\end{equation}   
A triple $(g,f,\tau)$ satisfying (\ref{soliton}) is called a shrinking Ricci soliton, the associated family $g(t)$ of metrics along (\ref{system}) 
is actually a self-similar solution\footnote{that is, a family of metrics moving by diffeomorphisms} of the Ricci flow.

\quad  We would now like to argue that the $\mathcal{W}$-entropy cannot arise from the scheme described in this note - at least not 
not by following the scheme exactly.  According to (\ref{dS}), the vanishing of $\frac{\partial S}{\partial \tau}$ happens 
if and only if all microstates share the same energy value.  On the other hand, if the Dirichlet energy 
$\langle\varphi, A_{\theta}\varphi\rangle$ is constant for all $\varphi$, then the constant must be zero.  For a self-adjoint operator $A_\theta$ 
this is impossible, since we immediately run into trouble for an eigenvector $\varphi$ of $A_\theta$.  Therefore, not only is it impossible 
to produce the $\mathcal{W}$-entropy using the functional determinant of self-adjoint operators as we have done, 
any entropy produced from our scheme can also never contain ``soliton-like'' objects, i.e. solutions $\theta(\tau)$ to 
the associated flow involved such that $\frac{dS}{d\tau} = 0$ along $\theta(\tau)$.  The formula  
(\ref{varS}) corroborates this claim except for $\chi(M) = 6$, although we do not know how one can construct a closed, connected surface with Euler 
Characteritic $6$.    

\quad In our opinion, the scheme described in this paper seem to be the most direct approach to constructing a partition function (and the ensuing 
entropy) associated to a geometric flow.  Therefore it is desirable 
that some modification of our scheme still works in constructing the partition function for the $\mathcal{W}$-entropy.  One possible modification 
is to consider functional determinants of operators that are not self-adjoint, but we do not know 
of any regularization schemes for such operators.\footnote{Note that a finite-dimensional example would be 
rotation by $90$ degrees about a fixed axis, whose Dirichlet energy is always zero since the rotated vector is orthogonal to the original vector.  
Thus macrostates that determine unitary operators could be ``solitons'' associated to some geometric flow, 
but due to their lack of real eigenvalues something like 
the zeta-function regularization seems elusive.}  Even if such an operator $A$ exists, there is still a problem of ensuring that it be 
invariant under the flow (\ref{system}).  In view of the cases involving conformal deformations that we have dealt with, it is also unclear to us whether or 
not to consider (\ref{system}) as part of a more general family of flows.\footnote{From all the work in the literature it seems that 
(\ref{system}), at least for a general manifold dimension $n$, is quite special in its own right.}   
Fitting the $\mathcal{W}$-entropy into our scheme seems to be a difficult problem. 

\quad One last remark before we leave this section.  In statistical mechanics, there is also the concept of the \text{\it{free energy}}, defined by  
 \[
  F = -\tau \log Z
 \]
with $\tau$ the temperature and $Z$ the partition function as before.  Clearly, the entropy (\ref{eqgibbs}) can be expressed using the free energy as 
\begin{equation}\label{SF}
 S = -\frac{\partial F}{\partial\tau}.
\end{equation}
Therefore, Perelman's declaration of (\ref{lnZ}) can be seen as a declaration of the free energy (up to a factor of $-\tau$).  A straight-forward exercise 
also confirms that the $\mathcal{W}$-entropy also follows from (\ref{SF}) by differentiating along the flow (\ref{system}).\footnote{Some 
mathematicians have already observed this mathematical fact, see for example \cite{KL}.}   

\section{Physical Interpretation}
\quad So far we have not explicitly discussed the thermodynamic system underlying our scheme.  Since the formulas 
we have conjured up are quite consistent with the statistical mechanics (in particular, the fact that (\ref{varS}) does not change sign), 
a more detailed physical interpretation should be in order.     
Recall the
general formula 
\[
 Z = \int_{\mathcal{H}} e^{-\beta E(\varphi)}\,\mathcal{D}\varphi
\]
of the partition function which we have employed throughout. 
In the most general sense, the microstates $\varphi$ should be interpreted as 
\text{\it{scalar fields}} and the Dirichlet energy $E(\varphi)$ can 
then be interpreted as the corresponding field energy.  $\mathcal{H}$ is then just the field configuration space.  This is the most abstract 
interpretation, as it covers all potential examples of operators that one can apply our scheme to.  On the other hand, the 
energy-determining macrostates are external (or background) constraints.  These energy-determining macrostates 
should not be thought of as being analogous to the volume of a system of gas particles, as there the volume affects both the nature and the 
number of microstates (in the quantum and classical description, respectively).  Since operators in our scheme for different energy-determining
macrostates share a dense intersection of their domains, the microstates can be fixed as this dense intersection.  Therefore, the best analogy for the 
energy-determining macrostates in our scheme would be to the external magnetic field applied to a system of a fixed number of 
particles with spin.  The spin states of the particles do not change with the strength of the magnetic field applied, but the ensuing energy does.

\quad In the special case of the Laplacian on closed surfaces, we speculate that a more concrete physical interpretation may be available.  
In this case, the microstates $\varphi$ can be thought of as the vibrational profile of a membrane on the surface.  
In particular, the value of $|\varphi|$ at a point on the surface represents the transverse displacement of the membrane from the surface.  
The transverse displacement of the membrane should correspond to some kind of ``stretching'' of the membrane,\footnote{Note that this is 
different from the stretching that one considers for a general \text{\it{harmonic map}}, where the stretching of the membrane is 
along the intrinsic, tangent directions of the manifold.}
 and the quantity $|\nabla\varphi|^2$ 
should be a measure of the magnitude of such a stretching - as seen by the Riemannian metric on the surface.  Then the Dirichlet energy 
can be seen as the total amount of energy stored as a result of such a stretching.  This however, only takes care of the potential 
energy.  The correct physical interpretation should also include a kinetic part of the energy, which is not included in the Dirichlet energy.

\section{Appendix}
\quad By the finite-dimensional version of (\ref{classic}) we mean pretending that the function space $\mathcal{H}$ is just $\mathbb{R}^n$ for some
finite $n$.  Suppose $A_g \psi = \lambda \psi$, then we must have $A_{\beta^{-1}g}\psi = \beta\lambda\psi$.  It follows 
immediately that a normalized eigenvector 
for the operator $A_{\beta^{-1}g}$ in the metric $\beta^{-1}g$ is $\tilde{\psi}_i = \sqrt{\beta}\psi$, where $\psi$ is a normalized eigenvector 
for $A_g$ in the metric $g$.  Then for a given vector $\varphi$, we have the expansion
\[
 \varphi = \tilde{c}_i\tilde{\psi}_i = \tilde{c}_i \sqrt{\beta}\psi_i = c_i \psi_i,
\]            
which implies that $\tilde{c}_i = \beta^{-1/2}c_i$.  

\quad Expanding $\varphi$ in $\psi_i$, we compute that $\langle\varphi, A_{\beta^{-1}g}\varphi\rangle = \beta\lambda_j c_j^2$ and the 
left-hand side of (\ref{classic}) becomes
\begin{align}
 \int_{\mathcal{H}} e^{-\langle\varphi, A_{\beta^{-1}g}\varphi\varphi\rangle_g}\, \mathcal{D}\varphi_g  &= 
 \int_{\mathbb{R}^n} e^{-\beta\lambda_j c_j^2}\,\pi^{-n/2}dc_1\cdots dc_n \notag\\
 &= \pi^{-n/2} \int_{\mathbb{R}} e^{-\beta\lambda_1 c_1^2}\,dc_1 \cdots \int_{\mathbb{R}} e^{-\beta\lambda_n c_n^2} \,dc_n \notag\\
 &= \pi^{-n/2} \sqrt{\frac{\pi}{\beta\lambda_1}}\cdots\sqrt{\frac{\pi}{\beta\lambda_n}} = \frac{\beta^{-n/2}}{\sqrt{\lambda_1\cdots\lambda_n}}.\notag
\end{align}
If we expand $\varphi$ in $\tilde{\psi}_i$ and use the formula $\mathcal{D}\varphi_{\beta^{-1}g} = \pi^{-n/2}d\tilde{c}_i\cdots d\tilde{c}_n$ with 
$d\tilde{c}_i = \beta^{-1/2} dc_i$, then the righ-hand side of (\ref{classic}) becomes 
\begin{align}
 \int_{\mathcal{H}} e^{-\langle\varphi, A_g \varphi\rangle_g}\,\mathcal{D}\varphi_{\beta^{-1}g} &= 
 \int_{\mathbb{R}^n} e^{-\lambda_i c_i^2} \, \pi^{-n/2}\beta^{-n/2} dc_1\cdots dc_n \notag\\
 &= \pi^{-n/2}\beta^{-n/2}\int_{\mathbb{R}} e^{-\lambda_1 c_1^2}\,dc_1 \cdots \int_{\mathbb{R}}e^{-\lambda_n c_n^2}\,dc_n \notag\\
 &= \pi^{-n/2}\beta^{-n/2} \sqrt{\frac{\pi}{\lambda_1}}\cdots\sqrt{\frac{\pi}{\lambda_n}} = \frac{\beta^{-n/2}}{\sqrt{\lambda_1\cdots\lambda_n}}.\notag
\end{align}
Thus equation (\ref{classic}) holds, at least in finite dimensions.  In the calculations above, we have used the well-known integral
\[
 \int_{-\infty}^{\infty} e^{-ax^2}\,dx = \sqrt{\frac{\pi}{a}}\,.
\]
Note that the coordinates $(c_1,...,c_n)$ coming from eigenvectors of $A_g$ are global on the Euclidean space $\mathbb{R}^n$ (viewed as a smooth manifold), 
and therefore $\pi^{-n/2} dc_1\cdots dc_n$ are different measures on $\mathbb{R}^n$ for different choices of $g$.  This explains why we should a-priori 
assume that the infinite-dimensional version of the measure is also dependent on the energy-determining macrostates.     

\quad Next we examine the finite-dimensional validity of $J(g_0, g)$ appearing in formula (\ref{additional}).  
The argument below generalizes the one above for conformally scaling a metric by a constant.  As we did above, we are taking 
$\mathcal{H} = \mathbb{R}^n$.  Let $g$ and $g_0$ be two different metrics, and for a given $\varphi\in\mathcal{H}$ we can expand it as 
$\varphi = c_i\psi_i$ with $\psi$ the normalized eigenvectors according to $g_0$, or expand it as $\varphi = \tilde{c}_i\tilde{\psi}_i$ with 
$\tilde{\psi}_i$ the normalized eigenvectors according to $g$.  Then we can write $\tilde{\psi}_i = a_{ij}\psi_j$ with the matrix $\big(a_{ij}\big)$ 
with $\big(a^{ij}\big)$ its inverse, which implies $\tilde{c}_i = a^{ij} c_j$ and hence $d\tilde{c}_i = a^{ij} dc_j$.  From this we see that  
\begin{align}
 \mathcal{D}\varphi_g &= \pi^{-n/2}d\tilde{c}_1\cdots d\tilde{c}_n \notag\\
  &= \pi^{-n/2} \det\big(a^{ij}\big) dc_1\cdots dc_n = \det\big(a^{ij}\big)\, \mathcal{D}\varphi_{g_0}.
\end{align}   
Therefore (\ref{additional}) is valid in the finite-dimensional case, with $J(g_0,g) = \det\big(a^{ij}\big)$.  Going from 
$\{\tilde{\psi}_i\}$ to $\{\psi_i\}$ is a linear coordinate-change on $\mathcal{H}$, and thus $\det\big(a^{ij}\big)$ can be interpreted as 
the Jacobian of the coordinate-change.  In particular,   
$\det\big(a^{ij}\big)$ is a function that depends only on $g_0$ and $g$ (and not on $\varphi$).

\vskip 0.8cm

\text{\it{E-mail Address}}: ccl37@case.edu

\end{document}